\documentclass[12pt]{article}

\font\tenmsb=msbm10
\font\sevenmsb=msbm7
\font\fivemsb=msbm5

\newfam\msbfam
\textfont\msbfam=\tenmsb
\scriptfont\msbfam=\sevenmsb
\scriptscriptfont\msbfam=\fivemsb
\def\Bbb#1{{\fam\msbfam #1}}

\makeatletter
\ifnum\@ptsize=0  \fi \ifnum\@ptsize=2
 \fi 

\newcommand\bR{{\Bbb R}}
\newcommand\bZ{{\Bbb Z}}
\newcommand\bC{{\Bbb C}}
\newcommand\bQ{{\Bbb Q}}

\newcommand\bP{{\Bbb P}}

\newtheorem{theorem}{Th\' eor\`eme}[section]
\newtheorem{lemma}[theorem]{Lemme}
\newtheorem{corollary}[theorem]{Corollaire}
\newtheorem{proposition}[theorem]{Proposition}

\newtheorem{re}[theorem]{Remarque}

\newtheorem{definition}[theorem]{Definition}

\newenvironment{remark}{\begin{re}\em}{\end{re}}

\begin{document}
\title {VARI\'ET\'ES FAIBLEMENT SP\'ECIALES \`A COURBES ENTI\`ERES D\'EG\'EN\'ER\'EES} \author{Fr\'ed\'eric Campana et Mihai P\u aun}

\maketitle

*
{\bf Summary:} A complex projective manifold is said to be {\it weakly-special} if none of its finite \'etale covers has a dominat rational map to a positive-dimensional manifold of general type. Weakly-special manifolds are conjectured in [H-T 00] to be potentially dense if defined over a number field. In [Ca 01; 2.1], the stronger notion of special (complex projective) manifold is introduced. These are conjectured there (9.20) to be {\it exactly} the potentially dense ones (when defined over a number field). The two notions of specialness and weak-specialness coincide up to dimension $2$, but differ from dimension $3$ on, as shown by examples $X$ constructed by Bogomolov-Tschinkel in [BT 03], answering a question raised in [Ca 01]. These examples should thus be potentially dense, according to the first conjecture, but should not be so, according to the second conjecture. The present techniques of arithmetic geometry do not seem to be able to decide between these two conflicting conjectures. So, instead of this, we investigate the more tractable hyperbolic aspects of the examples of [B-T 03], and show that for certain $X$'s, the behaviour of entire transcendental curves is the one conjectured in [Ca 01; 9.16, p. 618]. More precisely: the image of any holomorphic map from $\bC$ to $X$ is contained either in some fibre of the unique elliptic fibration on $X$, or in some fixed divisor of $X$.  Which is consistent with the standard links (Lang, Vojta for varieties of general type, [Ca 01] in general) expected between hyperbolic and arithmetic behaviour of projective varieties in case the second conjecture (in [Ca 01]) holds. Conversely,  if the conjecture stated in [H-T 00] is true, these expected links between arithmetics and complex hyperbolicity should fail to hold on these examples already.

\


\


\

\

\vspace*{-0.5in}\section*{Introduction}

Soit $X$ une vari\'et\'e projective complexe lisse et connexe. Une fibration $\varphi:X\to Y$ est une application m\'eromorphe dominante \`a fibres connexes, avec $Y$ lisse.

On dit que $X$ est {\it faiblement-sp\'eciale} (ou {\it f-sp\'eciale}) s'il n'existe pas de rev\^etement \'etale fini $u:X'\to X$ et de fibration $\varphi:X'\to Y$, o\`u $Y$ est une vari\'et\'e de type g\'en\'eral de dimension $p>0$ (voir [Ca 01; 9.26, p. 621] pour cette notion).

On dit que $X$ est {\it sp\'eciale} s'il n'existe pas de {\it fibration de type g\'en\'eral} $\varphi:X\to Y$, avec $dim(Y):=p>0$. Voir plus de d\'etails dans le \S1 ci-dessous.

La notion de fibration de type g\'en\'eral est d\'efinie dans [Ca01; 2.1,p.527]. Elle signifie que la base orbifolde de $\varphi$ d\'eduite de ses fibres multiples est de type g\'en\'eral pour tout mod\`ele bim\'eromorphe de $\varphi$. (Si $\varphi$ n'a pas de fibre multiple en codimension $1$, $\varphi$ est de type g\'en\'eral si et seulement si $Y$ l'est). Voir \S 1 ci-dessous pour ces notions.

Une vari\'et\'e sp\'eciale est $f$-sp\'eciale (voir [Ca01; 1.8, 9.27]). R\'eciproquement, les courbes et surfaces sp\'eciales sont $f$-sp\'eciales. Des exemples de vari\'et\'es faiblement sp\'eciales non-sp\'eciales de dimension $3$ sont cependant construites par Bogomolov-Tschinkel dans [B-T03], r\'epondant ainsi affirmativement \`a une question pos\'ee dans [Ca01]. La construction de [B-T 03] est rappel\'ee dans \S1 ci-dessous. Les exemples de Bogomolov-Tschinkel sont simplement connexes et munis d'une unique fibration elliptique $\varphi:X\to B$, o\`u $B$ est une surface avec $\kappa(B)=1$; cette fibration $\varphi$ est cependant de type g\'en\'eral, car elle a des fibres multiples de multiplicit\'e $m\geq 2$ au-dessus d'une courbe $D$ telle que $\kappa(B,K_B+(1-1/m).D)=2$. Mais ces fibres multiples ne peuvent \^etre \'elimin\'ees par un rev\^etement \'etale de $X$, qui est simplement connexe.

Pour $X$ d\'efinie sur un corps de nombres $k$, on dispose de $2$ conjectures concernant la densit\'e potentielle de $X$ (ie: le fait que $X(K)$ soit Zariski dense dans $X$ pour une extension finie assez grande $K/k$):

A. Il est conjectur\'e dans [H-T00] (conjecture attribu\'ee \`a D.Abramovich et J.L. Colliot-Th\'el\`ene) que $X$ est potentiellement dense si et seulement si $X$ est faiblement sp\'eciale. Cette conjecture affirme donc que les $X$ de [B-T03] sont potentiellement denses (on peut en d\'efinir sur des corps de nombres). 

B. Il est ind\'ependamment conjectur\'e dans [Ca01; 9.20, p.619] que $X$ est potentiellement dense si et seulement si $X$ est sp\'eciale. Cette conjecture affirme donc que les $X$ de [B-T03] ne sont pas potentiellement denses

Les techniques actuelles de g\'eom\'etrie arithm\'etique ne semblent pas pouvoir d\'ecider entre ces deux affirmations contradictoires. On va donc consid\'erer ici les propri\'et\'es, plus accessibles, d'hyperbolicit\'e complexe des exemples de Bogomolov-Tschinkel. Ces propri\'et\'es sont conjecturalement (Lang, Vojta pour les vari\'et\'es de type g\'en\'eral, [Ca 01] pour les vari\'et\'es arbitraires) de nature g\'eom\'etrique, et \'equivalentes aux propri\'et\'es qualitatives de distribution des points $K$-rationnels, apr\`es extension assez grande $K$ du corps de d\'efinition.

Nous allons montrer que, pour certains des exemples de [B-T 03], les courbes enti\`eres transcendantes $h:\bC\to X$ trac\'ees sur $X$ sont alg\'ebriquement d\'eg\'en\'er\'ees et que leur image est, soit contenue dans un diviseur $Z$ fix\'e de $X$, soit  contenue dans une fibre de $\varphi:X\to B$.

Si la conjecture A. est vraie, cette propri\'et\'e montre que les correspondances conjecturales attendues entre hyperbolicit\'e et arithm\'etique sont fausses. Par contre, cette propri\'et\'e d'hyperbolicit\'e est exactement celle conjectur\'ee dans [Ca 01] , et est compatible avec la conjecture B. et les correspondances attendues entre hyperbolicit\'e et arithm\'etique.

La conjecture $H$ de [Ca01; 9.2, p. 614] affirme en effet que $X$ est sp\'eciale si et seulement si sa pseudom\'etrique de Kobayashi est identiquement nulle, et dans le cas non-sp\'ecial des $X$ pr\'ec\'edentes, que la pseudom\'etrique de Kobayashi de $X$ est l'image r\'eciproque par $\varphi$ d'une pseudom\'etrique sur $B$, qui est une m\'etrique sur un ouvert de Zariski dense de $B$. (Ceci peut \^etre vu comme une version orbifolde de la Conjecture hyperbolique de S.Lang pour les surfaces de type g\'en\'eral). En particulier, pour les $X$ de Bogomolov-Tschinkel, la projection par $\varphi$ des courbes enti\`eres trac\'ees sur $X$ doit \^etre, soit r\'eduite \`a un point, soit contenue dans un ensemble alg\'ebrique strict de $B$, qui est donc une r\'eunion finie de courbes rationnelles ou elliptiques. C'est en effet ce que nous d\'emontrons.

La m\'ethode de d\'emonstration consiste \`a exploiter le fait que la base orbifolde $(B/\Delta)$ de $\varphi$ est de type g\'en\'eral, avec $\Delta:=(1-1/m).D$. Si le diviseur $D$ de $B$ \'etait divisible par $m$ dans $Pic(B)$, cette orbifolde serait rev\^etue de mani\`ere \'etale au sens orbifolde par un rev\^etement cyclique $\tau:B'\to B$ de degr\'e $m$ de $B$ ramifiant \`a l'ordre $m$ exactement le long de $D$, $B'$ \'etant une surface de type g\'en\'eral, puisque $K_{B'}=\tau^*(K_B+\Delta)$. Le changement de base normalis\'e $B'\to B$ fournirait alors $\varphi':X':=X\times _{B}B'\to B'$ , avec $X'$ \'etale sur $X$, donc tel que les courbes rationnelles ou elliptiques (resp. les courbes enti\`eres) trac\'ees sur $X$ se rel\`event \`a $X'$. La question se ram\`enerait donc, dans ce cas, apr\`es projection sur $B'$ par $\varphi$, \`a montrer la finitude des courbes rationnelles ou elliptiques sur $B'$, ainsi que la d\'eg\'en\'erescence alg\'ebrique des courbes enti\`eres trac\'ees sur $B'$. Si, de plus, le couple $(B,D)$ \'etait choisi tel que $(c_1^2-c_2)(B')>0$, on d\'eduirait de [B 77] (resp. [MQ 98]) la premi\`ere (resp. la seconde) propri\'et\'e. 

Cependant, le fait que $X$ soit $f$-sp\'eciale est exactement ce qui interdit \`a cette approche de fonctionner telle quelle (observer que nous pouvons contourner la difficult\'e dans le \S3 par des choix ad\'equats pour la premi\`ere propri\'et\'e, mais non pour la seconde). La propri\'et\'e cruciale utilis\'ee (et qui traduit l'existence de fibres multiples de $\varphi)$ est la suivante: si $h:\bC\to X$ est une courbe enti\`ere sur $X$, alors sa projection $\varphi\circ h:\bC\to B$ est $m$-tangente  \`a $D$, ce qui signifie qu'en chaque point d'intersection de l'image de cette courbe avec $D$, l'ordre de contact est divisible par $m$ (on pourrait se contenter de {\it au moins} $m$, qui est la bonne notion). On montre alors que les courbes enti\`eres de $B$ qui sont $m$-tangentes \`a $D$ sont alg\'ebriquement d\'eg\'en\'er\'ees si $(c_1^2-c_2)(B/\Delta)>0$(dans un sens ad\'equat). La m\'ethode consiste \`a d\'efinir et montrer l'existence de sections non nulles de $[Sym^N](\Omega_{(B/\Delta)}^1)\otimes A^{-1}$, $A$ ample sur $B$, pour $N$ assez grand, sous ces hypoth\`eses. Le symbole $[Sym^N]$ d\'esigne la {\it partie enti\`ere} du faisceau consid\'er\'e (qui n'a de sens que localement, une fois une racine $m$-i\`eme d'un \'equation locale de $D$ ayant \'et\'e choisie). Des arguments classiques de th\'eorie de Nevanlinna, compl\'et\'es par des r\'esultats r\'ecents concernant les courants d\'efinis par les courbes enti\`eres trac\'ees sur les surfaces ([MQ98],[Br00],[Pa03]) s'adaptent alors \`a notre contexte orbifolde et permettent de conclure. Voir  le \S 4 pour les d\'etails.

\

Nous remercions L. Caporaso pour nous avoir indiqu\'e de nombreuses imperfections dans la premi\`ere version de ce texte.

\

\newpage

\tableofcontents

\section{Vari\'et\'es sp\'eciales}

Nous rappelons bri\`evement les notions introduites dans [Ca 01], et qui jouent un r\^ole pour la compr\'ehension de la suite. 

Soit $X$ une vari\'et\'e projective complexe lisse et connexe. Une fibration (r\'eguli\`ere) $f:X\to Y$ est un morphisme de $X$ sur une vari\'et\'e projective lisse $Y$, \`a fibres connexes. Une modification de $f$ est une fibration $f':X'\to Y'$ avec $u:X'\to X,v:Y'\to Y$des modifications propres, et telles que $fu=v f'$.

La base orbifolde de $f$, not\'ee $\Delta(f):=\sum_D (1-1/m(f,D)).D$ est le $\bQ$-diviseur (effectif, fini) sur $Y$ d\'efini comme suit: $D$ varie parmi tous les diviseurs effectifs irr\'eductibles de $Y$, et pour un tel $D$, $m(f,D)\geq 1$ est l'entier d\'efini par $m(f,D):=inf \lbrace m_j,j\in J(D,f)\rbrace$.
On a utilis\'e les notations suivantes: $f^*(D)=\sum_{j\in J(D,f)}m_j.E_j+R$, o\`u $R$ est un diviseur de $X$ tel que $f(R)$ soit de codimension au moins $2$ dans $Y$, et les $E_j$ sont les diviseurs irr\'eductibles distincts de $X$ tels que $f(E_j)=D,\forall j$.

On d\'efinit alors: 

1. Le fibr\'e canonique de $(Y/\Delta(f)):=K_Y+\Delta(f)$: c'est un $\bQ$-diviseur sur $Y$.

2.  $\kappa(Y\Delta):=\kappa(Y,K_Y+\Delta)\geq \kappa(Y)$. 

Si $f':X'\to Y'$ est une modification de $f$, on a: $\kappa(Y'/\Delta(f'))\leq \kappa(Y/\Delta(f))$. On peut avoir in\'egalit\'e stricte si $\kappa(Y)=-\infty$ (mais pas si $\kappa(Y)\geq 0$). 

On dit que $f$ est une {\it  fibration de type g\'en\'eral } si $\kappa(Y'/\Delta(f'))=dim(Y)>0$ pour toute modification $f'$ de $f$. Lorsque $\kappa(Y)\geq 0$ (et c'est le cas des exemples de Bogomolov-Tschinkel), il suffit de consid\'erer $f$ seulement, la dimension canonique de la base orbifolde de $f'$ ne d\'ependant alors pas du mod\`ele bim\'eromorphe $f'$ de $f$ consid\'er\'e. Voir [Ca 01, proposition 1.14, p. 514). Cette notion est bim\'eromorphe, et permet de d\'efinir $\kappa(f):=\kappa(Y'/\Delta(f'))$ lorsque $f:X\to Y$ est maintenant une application m\'eromorphe dominante \`a fibre g\'en\'erique connexe, entre vari\'et\'es projectives irr\'eductibes \`a singularit\'es arbitraires, et $f':X'\to Y'$ un mod\`ele bim\'eromorphe de $f$ pour lequel la base orbifolde de $f'$ est de dimension canonique minimum.

\begin{definition}On dit que $X$ est {\it sp\'eciale } s'il n'existe pas de fibration m\'eromorphe de type g\'en\'eral $f:X\to Y$ 

\end{definition}

Les exemples fondamentaux de vari\'et\'es sp\'eciales sont les vari\'et\'es rationnellement connexes et les vari\'et\'es $X$ avec $\kappa(X)=0$. Inversement toute $X$ sp\'eciale peut \^etre canoniquement d\'ecompos\'ee en tour de fibrations  de l'un de ses deux types (dans un sens orbifolde ad\'equat, sous r\'eserve de la validit\'e de la conjecture $C_{n,m}^{orb}$).

On conjecture dans [Ca 01; 9.2, 9.20] que les vari\'et\'es sp\'eciales sont les analogues de dimension sup\'erieure des courbes rationnelles ou elliptiques:

1. $X$ est sp\'eciale si et seulement si sa pseudom\'etrique de Kobayashi est nulle.

2. $X$ est sp\'eciale si et seulement si elle est potentiellement dense (si d\'efinie sur un corps de nombres)

3. Si $X$ est sp\'eciale, $\pi_1(X)$ est virtuellement ab\'elien.

\

L'un des r\'esultats principaux de [Ca 01; 5.8, 5.16, pp. 579-582] est:

{\begin{theorem} Soit $X$ projective, irr\'eductible. Il existe une unique fibration $c_X:X\to C(X)$, appel\'ee le {\it coeur} de $X$, telle que:

1. Les fibres g\'en\'erales de $X$ sont sp\'eciales.

2. $c_X$ est soit une fibration de type g\'en\'eral , soit constante (ceci si et seulement si $X$ est sp\'eciale).

\end{theorem}

Le coeur scinde donc $X$ en ses deux composantes: sp\'eciale (les fibres) et de type g\'en\'eral (la base orbifolde). Il fournit une description synth\'etique de la g\'eom\'etrie de $X$ analogue \`a celle des surfaces due \`a Enriques-Shafarevitch-Kodaira.

Ce scindage a lieu aussi, conjecturalement, aux niveaux hyperbolique et arithm\'etique. On conjecture en effet dans [Ca 01] en particulier qu' il existe un sous-ensemble alg\'ebrique strict $Z\subset C(X)$ tel que:

1. Si $h:\bC\to X$ est une courbe enti\`ere de $X$, alors $c_X\circ h:\bC\to C(X)$ est soit constante, soit d'image contenue dans $Z$.

2. Si $K$ est un corps de nombres sur lequel $X$ est d\'efinie (ainsi donc que $c_X$), alors l'intersection de $c_X(X(K))$  avec $U:=X-Z$ est fini.

\

Ces conjectures sont inspir\'ees par celles de Lang pour les vari\'et\'es de type g\'en\'eral, et en sont la version orbifolde, par composition avec le coeur.

\newpage

\section{ La construction de Bogomolov-Tschinkel}

\

\subsection{Les exemples de Bogomolov-Tschinkel}

\begin{theorem}\label{bt} Il existe une fibration elliptique $\varphi:X\to B$ dans laquelle $X,B$ sont projectives complexes lisses et connexes, $dim(X)=3$, $dim(B)=2$, telle que de plus:

1. $\pi_1(X)=\lbrace 1\rbrace$

2. $\kappa(S)=1$

3. $\kappa(B,K_B+\Delta(\varphi))=2$.

Une telle  fibration peut \^etre d\'efinie sur $\Bbb Q$.

\end{theorem}

Une construction sera d\'ecrite dans la section suivante. On en d\'eduit:

\begin{corollary} Pour $\varphi:X\to  B$ comme ci-dessus: $X$ est f-sp\'eciale, non sp\'eciale, et $\varphi=c_X$ est le coeur de $X$.
\end{corollary}

{\bf D\'emonstration:} La fibration $\varphi$ est \`a fibres g\'en\'eriques sp\'eciales (courbes elliptiques), et base orbifolde de type g\'en\'eral (par la propri\'et\'e 3.). C'est donc le coeur de $X$, qui n'est donc pas sp\'eciale. Pour montrer que $X$ est faiblement sp\'eciale, il suffit, puisqu'elle est simplement connexe, de voir qu'elle n'a pas de fibration $f:X\to Y$ sur une base $Y$ de type g\'en\'eral avec $d:=dim(Y)>0$. On ne peut avoir ni $d=3$ (puisque $X$ est couverte par des courbes elliptiques), ni $d=1$ (puisque $X$ est simplement connexe). Si $d=2$, on a $f\neq \varphi$, puisque $\kappa(B)=1$; donc $Y$ est couverte par les images des fibres de $\varphi$, et ne peut donc \^etre de type g\'en\'eral. Contradiction.

\

L'objectif du pr\'esent texte est de montrer le:

\

\begin{theorem} Il existe des fibrations $\varphi:X\to B$ poss\'edant les propri\'et\'es \'enonc\'ees en \ref{bt}, telles qu'existe une courbe projective $\Gamma\subset B$ telle que pour toute courbe enti\`ere $h:\bC\to X$, l'image de $\varphi\circ h:\bC\to B$ est soit un point, soit contenue dans $\Gamma$ (``hyperbolicit\'e transcendante"). 

En particulier, si $C\subset X$ est une courbe rationnelle ou elliptique trac\'ee sur $X$, alors $\varphi(C)\subset B$ est soit un point, soit contenue dans $\Gamma$ (``hyperbolicit\'e alg\'ebrique").

On peut trouver de tels $\varphi$ d\'efinies sur $\bQ$.

\end{theorem}

L'hyperbolicit\'e alg\'ebrique sera \'etablie dans la \S \ref{halg} par un argument tr\`es simple. La version transcendante le sera dans \S\ref{htransc}.

\

\subsection{ La Construction}\label{bt}

Nous conservons les notations de [B-T03]. Soit $(B,\varphi_1,\psi,\eta',E')$ un quintuplet tel que:

\

BT1. $\eta':E'\to C$ est une surface elliptique projective ayant une (unique) fibre multiple $(\eta')^{-1}(0):=E'_0$, de multiplicit\'e $m\geq 2$, fix\'ee dans la suite du texte. On suppose aussi que $E'-E'_0$ est simplement connexe. On peut obtenir de telles donn\'ees par transformation logarithmique d'une surface rationnelle elliptique ad\'equate.

\

BT2. Une surface elliptique minimale projective $\varphi_1:B\to C_1\cong \bP^1$ avec $\kappa(S)=1$, telle que, de plus:

\

BT3. il existe une fibration $\psi:B\to \bP_1$, distincte de $\varphi_1$, dont une fibre lisse $D:=\psi^{-1}(0)$ (n\'ecessairement de genre $g\geq 2)$ soit lisse et telle que $B-D$ soit simplement connexe.

\

La construction de [B-T03] consiste alors \`a prendre pour $\varphi:X:=E'\times_{\bP_1}B\to B$ la fibration d\'eduite de  $\eta':E'\to \bP_1$ par le changement de base $\psi:B\to \bP_1$, apr\`es normalisation du produit fibr\'e. 

On a donc: $\Delta(\varphi)=\psi^{-1}(\Delta(\eta'))=\psi^{-1}((1-1/m).E'_0)=(1-1/m)\varphi^{-1}(D)$. On en d\'eduit imm\'ediatement les propri\'et\'es \'enonc\'ees en \ref{bt} (la simple connexit\'e de $X$ r\'esulte de celle de $X-\varphi^{-1}(D))$.

\

Montrons maintenant, suivant [B-T03], comment r\'ealiser la simple-connexit\'e de $B-D$. Dans la suite, nous consid\'ererons toujours cette situation particuli\`ere.

\

BT4. Soit $B_0$ une surface projective lisse et simplement connexe telle que $\kappa(B_0)=1$. Soit $D_0,D_0'\in \vert L_0\vert $ deux membres g\'en\'eriques (et m\^eme \'eventuellement g\'en\'eraux), o\`u $L_0\in Pic(B_0)$ est ample, et $\vert L_0\vert$ soit sans point base.

\

Soit alors $\psi:B_0\to C\cong \bP^1$ l'application m\'eromorphe d\'efinie par le syst\`eme lin\'eaire de dimension $1$ engendr\'e par $D_0$ et $D_0'$, qui se coupent transversalement en un ensemble fini not\'e $F$. On note $\beta:B\to B_0$ l'\'eclatement de $B_0$ le long de $F$, et $D\subset B$ le transform\'e strict de $D_0$.

\

Montrons la simple connexit\'e de $(B-D)$. Par le th\'eor\`eme de Lefschetz, $\pi_1(B_0-D_0)\cong \bZ_r$, si $r>0$ est le plus grand entier par lequel $L_0$ est divisible dans $Pic(B_0)$. Et $\pi_1(B_0-D_0)$ est engendr\'e par un lacet $g$ tournant une fois autour de $D_0$. Le relev\'e par $\beta$ de ce lacet engendre donc $\pi_1(B-D)$, mais devient homotope \`a z\'ero dans $P-\lbrace u\rbrace\cong \bC$, si $P$ est une courbe exceptionnelle de $\beta$, et $u$ son point d'intersection avec $D$. D'o\`u la simple connexit\'e.

\

On fixera toujours $S_0,L_0,D_0$ dans la suite. Par contre, $D'_0$ sera variable et devra \^etre choisi g\'en\'erique (ou g\'en\'eral) dans le syst\`eme lin\'eaire d\'efini par $L_0$. On notera donc quelquefois $B_F$ la surface $B$, et $\varphi_F:X_F\to B_F$ la fibration associ\'ee pour indiquer la d\'ependance en $F$.

\

\begin{remark} Dans [B-T 03], les auteurs imposent au pinceau d\'efini par $D_0,D_0'$ d'\^etre tr\`es ample, et d'avoir tous ses membres irr\'eductibles. La seconde condition n'est pas n\'ecessaire. La premi\`ere peut \^etre affaiblie en: $D_0,D_0'$ sont amples et membres g\'en\'eriques d'un syst\`eme lin\'eaire sans point base (comme ci-dessus). La seule propri\'et\'e utilis\'ee est en effet que le groupe fondamental du compl\'ementaire de $D_0$ dans $S$ soit fini et engendr\'e par un petit lacet tournant une fois autour de $D_0$, \'egalement satisfaite sous cette hypoth\`ese plus faible.
\end{remark}

\

Dans la suite, nous aurons besoin de propri\'et\'es additionnelles BT 5 et BT6, que nous introduisons maintenant. 

\

BT5. On supposera toujours que $L_0$ est divisible par $m\geq 2$ (la multiplicit\'e de $E'_0$ dans $\eta')$. Il existe donc un rev\^etement cyclique $r:B'\to B_0$, de degr\'e $m$, ramifi\'e au-dessus de $D_0$ exactement. Puisque $D_0$ est lisse, $B'$ l'est aussi.

\

BT6. On supposera aussi que $(c_1^2-c_2)(B')>0$.

\

Montrons comment r\'ealiser BT6 aussi.

\

\subsection{ Exemples}\label{ex}

 Nous appliquerons deux fois le lemme suivant:

\begin{lemma} Soit $S$ une surface complexe lisse, compacte et connexe. Soit $D\in \vert mL\vert$ un diviseur lisse, et $r:S'\to S$ le rev\^etement cyclique de degr\'e $m$ de $S$ ramifi\'e le long de $D$. Les nombres de Chern $c_1^2(S')$ et $c_2(S')$ de $S'$ satisfont alors les relations suivantes (ceux de $S$ \'etant not\'es $c_1^2$ et $c_2$):

1. $c_1^2(S')=m. [c_1^2+(m-1).(2.K_S.L+(m-1). L^2]$.

2. $c_2(S')=m. [c_2+(m-1).(K_S.L+m. L^2)]$.

3. $c_1^2(S')-c_2(S')=m. [(c_1^2-c_2)+(m-1)(K_S.L- L^2)]$.

\end{lemma}\label{lemm}

{\bf D\'emonstration:} On a, par formule de ramification:

$K_{S'}=r^*(K_S+(m-1).L$, donc $c_1^2(S')=m.[c_1^2+2(m-1).K_S.L+(m-1)^2L^2]$.

Pour calculer $c_2(S')$, on remarque que $\chi_t=c_2$ \'etant la caract\'eristique topologique, la suite exacte de Mayer-Vietoris fournit: $c_2(S')=\chi_t(S')=m.\chi_t(S)-(m-1).\chi_t(D)=m.c_2(S)+(m-1).[K_S.D+D^2]=m.c_2(S)+(m-1).m.[K_S.L+m.L^2]$.

\

\begin{corollary} Dans la situation du lemme pr\'ec\'edent, $c_1^2(S')-c_2(S')>0$ si les conditions suivantes P1 et P2 sont satisfaites:

P1.  $K_S.L>L^2$

P2. $m>(1-[(c_1^2-c_2)/(K_S.L-L^2)])$.
\end{corollary}

Nous allons montrer sur des exemples tr\`es simples que les conditions du corollaire pr\'ec\'edent peuvent \^etre r\'ealis\'ees avec des surfaces elliptiques minimales dans la situation de \ref{bt}.

\

Soit $s:B_0\to P:=\bP^1\times\bP^1$ le rev\^etement double ramifi\'e le long du diviseur lisse 
$R\in \vert \cal O$$_P(2(k+2),4)\vert $. 
On suppose $R$ (et donc $S$) lisse. Par le lemme 1.1, on  a donc: $c_1^2(B_0)=0$, et $c_2(B_0)=12\chi(\cal O$$_S)=12.(k+2)$.

On v\'erifie de plus que $B_0$ est simplement connexe, par [F-M98,II,2.2, Prop. 2.1]. 

On choisit ensuite $L_0:=s^*(L')=s^*(O_P(b,ba))$, pour $a,b>0$ entiers. On a donc: $K_{B_0}.L_0=2kba$, et $L_0^2=4b^2a$, de sorte que la propri\'et\'e P1 (resp. P2) pr\'ec\'edente est satisfaite si $k\geq (2b+1)$, et si $(m-1).a>12.((k+2)/(k-2b))$.

\section{Hyperbolicit\'e Alg\'ebrique}\label{halg}

\

{\bf Hypoth\`eses 3.0.} On fixe $\varphi:X\to B$, construite comme en \ref{bt}. On suppose que les conditions BT1-BT6 sont satisfaites. En particulier: $(c_1^2-c_2)(B')>0$. 

\

On appelle courbe rationnelle (resp. elliptique) une courbe projective irreductible complexe $C$ de genre g\'eom\'etrique $g(\hat C):=0$ (resp. $1$), si $\hat C$ est la normalis\'ee de $C$.

\

On note $R\subset B_0$ la r\'eunion (d\'enombrable) des courbes rationnelles ou elliptiques de $B_0$. Donc $R\cap D_0$ est un ensemble d\'enombrable puisque $g(D_0)>1$. Puisque $ L_0$ est, par hypoth\`ese, sans point base, lorsque $D'_0$ est g\'en\'eral dans $\vert L_0\vert$, $F:=D'_0\cap D_0$ est disjoint de $R\cap D_0$.

\begin{theorem} On suppose les hypoth\`eses 3.0 ci-dessus r\'ealis\'ees, et $D'_0$ disjoint de $(D_0\cap R)$.  Il existe une courbe projective $Z\subset B$ telle que, pour toute courbe rationnelle ou elliptique $C$ de $X$, $\varphi(C)\subset B$ est soit un point, soit contenue dans $Z$.
\end{theorem}

{\bf D\'emonstration:}

Soit $\varphi':X':=X\times_{S}S'\to S'$ la (normalis\'ee de la) fibration d\'eduite de $\gamma:=\beta\circ\varphi:X\to B_0$ par le changement de base $r:B'\to B_0$ d\'efini en \S 2.2, propri\'et\'e BT5. Soit $u:X'\to X$ la projection naturelle. Alors $u$ est \'etale au-dessus de $(B_0-F)$, si $F\subset B_0$ est l'ensemble fini intersection de $D_0$ et $D_0'$.

Par le choix de $D'_0$,  $\varphi(C)$ ne rencontre pas $F$, si $C$ n'est pas contenue dans une fibre de $\varphi$ (ce que nous supposerons d\'esormais).

Donc, $u$ est \'etale au-dessus de $C$, et $C':=u^{-1}(C)$ est une r\'eunion (finie) de courbes rationnelles ou elliptiques de $X'$. Il en est de m\^eme de l'image $C":=\varphi'(C')\subset S'$ de $C'$ par $\varphi':X'\to S'$. De plus, aucune composante irr\'eductible de $C"$ n'est un point de $B'$, puisque $C$ n'est pas contenue dans une fibre de $\varphi$.

Or, $c_1^2(S')-c_2(B')>0$, d'apr\`es \ref{lemm}. On d\'eduit de [B77] (voir aussi [D77]) qu'il existe une courbe $Z'\subset B'$, r\'eunion finie de courbes rationnelles ou elliptiques, qui contient $C"$. Donc $\varphi(C)\subset Z$, si $Z:=\beta^{-1}(r(Z'))$. 

\

\section{Hyperbolicit\'e, version transcendante}\label{htransc}

\subsection{Notion de $m$-tangence}

\

Soit $S_0,L_0, F,B,B',D,X,\varphi$ d\'efinis comme en \S 2.2. On suppose les propri\'et\'es BT1-BT6 satisfaites.

\

\begin{definition}
Soit $h:M\to X$ une application holomorphe d'une courbe complexe (vari\'et\'e analytique complexe lisse et connexe de dimension pure $1)$ dans $B$. On dit que $h$ (ou $M$, par abus de language) est  {\it (classiquement) $m$-tangente} \`a $D$ si l'ordre de contact de $h(M)$ et de $D$ en chacun de leurs points d'intersection est {\bf divisible} par $m$ (ceci apr\`es composition avec la normalisation de $M$). Plus pr\'ecis\'ement: si $h^*(D)$ est un diviseur effectif de $M$ dont les coefficients sont tous divisibles par $m)$. Si $C$ est une courbe complexe (\'eventuellement singuli\`ere) trac\'ee sur $X$, $h$ sera la compos\'ee de l'inclusion naturelle de $C$ dans $X$ avec une normalisation de $C)$.
\end{definition}

\

{\bf Remarque:} La propri\'et\'e de $m$-tangence non classique est d\'efinie de mani\`ere similaire, mais en \' exigeant seulement un ordre de contact {\bf au moins \'egal} \`a $m$. Bien que plus naturelle, elle n'est pas utilis\'ee ici. Voir [Ca 04] pour les diff\'erences (parfois essentielles) entre les deux notions.

\

L'observation triviale suivante permet de r\'eduire \`a l'orbifolde $(B/\Delta(\varphi))$ l'\'etude des courbes trac\'ees sur $X$.

\

\begin{theorem} On suppose les hypoth\`eses 3.0 ci-dessus r\'ealis\'ees. Soit $h:M\to X$ une application holomorphe  de la courbe complexe $M$ vers $X$, telle que $\varphi\circ h:M\to B$ soit non constante. Alors: $\varphi\circ h:M\to B$ est (classiquement) $m$-tangente \`a $D$.
\end{theorem}

{\bf D\'emonstration:} Soit $H:=\varphi^{-1}(D)$ l'image r\'eciproque r\'eduite de $D$. C'est un diviseur irr\'eductible, et $\varphi^*(D)=m.H$, puisque les fibres de $\varphi$ au-dessus des points de $D$ sont toutes irr\'eductibles de multiplicit\'e $m$. On a donc: $(\varphi\circ h)^*(D)=h^*(\varphi ^*(D))=h^*(m.H)=m.h^*(H)$, et l'assertion.

\

\subsection{Objectif} 

\

\begin{theorem} \label{hyptransc} Soit $S_0,L_0,D_0, B'$ poss\'edant les propri\'et\'es BT1-BT6 \'enonc\'ees en \S 2.2. Pour $D'_0$ g\'en\'erique dans $\vert L_0\vert$, soit $\varphi_F:X_F\to B_F$ (voir notations BT5 de la \S 2.2)) la fibration associ\'ee.

Alors il existe une courbe alg\'ebrique $\Gamma\subset B=B_F$ telle que pour toute courbe enti\`ere $h:\bC\to B$ $m$-tangente \`a $D$ est soit constante, soit d'image contenue dans $\Gamma$. 

En particulier, si $h:\bC\to X$ est une courbe enti\`ere, alors $\varphi\circ h:\bC\to B$ est soit constante, soit d'image contenue dans $\Gamma$.
\end{theorem}

\

Nous allons pour ceci adapter au cadre orbifolde des arguments de th\'eorie de Nevanlinna. Le point cl\'e (tout comme dans la section pr\'ec\'edente) est que l'orbifolde $(B/\Delta)$, de dimension $2$, avec $\Delta:=(1-1/m).D)$ est de type g\'en\'eral avec $(c_1^2-c_2)>0$ (dans un sens orbifolde ad\'equat).

Nous ne d\'efinirons ici que les notions de formes diff\'erentielles orbifolde effectivement utilis\'ees (bien qu'il soit possible de travailler de fa\c con similaire dans un contexte g\'en\'eral).

Expliquons pourquoi il n'est pas possible ici d'appliquer au rev\^etement cyclique ramifi\'e de $S$ le long de $D$ les arguments tr\`es simples de la section pr\'ec\'edente pour obtenir la version transcendante de l'hyperbolicit\'e: c'est que courbes enti\`eres de $B$ (en g\'en\'eral non alg\'ebriquement d\'eg\'en\'er\'ees) peuvent recouvrir $B$, et que l'on ne peut alors pas choisir $F=D\cap D'$ disjoint de toutes ces courbes enti\`eres. Et l'on ne peut alors pas n\'ecessairement relever \`a $X'$ (d\'eduit par changement de base par ce rev\^etement) les courbes enti\`eres de $X$.

On montre n\'eammoins l'\'existence d'un nombre fini de points de $D$ tels que toute courbe enti\`ere trac\'ee sur $B$ et $m$-tangente \`a $D$ coupe $D$ en ces points seulement au plus. Cette \'etape repose sur la th\'eorie de Nevanlinna, appliqu\'ee aux formes pluridiff\'erentielles orbifolde sur $(B/\Delta)$. On utilise, suivant l'argument original de Bogomolov, l'in\'egalit\'e $(c_1^2-c_2)>0$ (dans un sens orbifolde ad\'equat) pour produire des feuilletages sur $B$ auquel toute courbe enti\`ere (ou bien rationnel ou elliptique) doit \^etre tangente.

\

\subsection{Formes \`a p\^oles logarithmiques orbifoldes}

\

On se donne une surface projective complexe lisse et connexe $B$, une courbe lisse et irr\'eductible $D$ sur $B$, et une multiplicit\'e enti\`ere $m\geq 2$. Le diviseur orbifolde consid\'er\'e sur $B$ est $\Delta:=(1-1/m).D$.

Notre objectif est le corollaire 4.4 ci-dessous, qui fournit des conditions suffisantes pour l'existence de formes pluridiff\'erentielles sur l'orbifolde $(B/\Delta)$. 

\

\begin{corollary} Soit $B,D,m$ comme ci-dessus, et $A$ ``big"  sur $B$. On suppose que $\kappa(B)=1$, que $g(D)>1$, que $\cal O$$_D(D))\cong \cal O$$_D$, et que $(c_1^2-c_2)(B)+(1-1/m).K_B.D>0$. 

Alors, il existe $N>0$ tel que $H^0(B, [Sym^N](\Omega ^1_{(B/\Delta)})\otimes A^{-1})\neq 0$.
\end{corollary}

\

Remarquons que les hypoth\`eses du pr\'ec\'edent corollaire sont satisfaites lorsque les donn\'ees $B,D,m$ sont d\'eduites par \'eclatement de $F$  de donn\'ees $B_0,L_0,D_0, F$ satisfaisant les hypth\`eses BT1-BT6 de la \S 2.2.

\

Expliquons bri\`evement l'origine de la notion de forme pluridiff\'erentielle sur l'orbifolde $(B/\Delta)$ (sous la forme implicite rudimentaire utilis\'ee ici).

\

Si $c_1(D)$ est divisible par $m$ dans $NS(B)$, on peut consid\'erer le rev\^etement cyclique $r:B'\to B$ de degr\'e $m$ ramifi\'e le long de $D$. 

Si $(x,y)$ sont des coordonn\'ees locales holomorphes pr\`es d'un point $a\in D$, en lequel $D$ a pour \'equation $y=0$, et si $a':=r^{-1}(a)\in S'$, alors les sections locales de $Sym^N(\Omega_{B'}^1)$ s'\'ecrivent sous la forme $$\sum_{j=0}^{j=N}a_j.(d(y^{1/m}))^{\otimes j}\otimes dx^{\otimes (N-j)}=\sum_{j=0}^{j=N}a_j.(dy/y^{(1-1/m)})^{\otimes j}\otimes dx^{\otimes (N-j)},$$ o\`u $y^{1/m}$ est une coordonn\'ee locale de $B'$ en $a'$ telle que $(y^{1/m})^m=r^*(y)$.

\

Ceci sugg\`ere de prendre la {\it partie enti\`ere} de cette expression, c'est-\`a-dire la somme des termes pr\'ec\'edents, mais en y rempla\c cant les exposants fractionnaires $(1-1/m).j$ par leur partie enti\`ere $[j-j/m]$. On obtient ainsi par image r\'eciproque par $r$ des sections de $Sym^N(\Omega_{B'}^1)$.

\

\begin{definition} Soit $[Sym^N](\Omega_{(B/\Delta)}^1)$ le faisceau des germes de sections m\'eromorphes de $Sym^N(\Omega_B^1)$ qui sont holomorphes en dehors de $D$, et qui s'\'ecrivent en coordonn\'ees locales $(x,y)$ pr\`es de $a\in D$ sous la forme:

$$\sum_{j=0}^{j=N}a_j.(dy^{\otimes j}/y^{[(1-1/m)j]})\otimes dx^{\otimes (N-j)}=\sum_{j=0}^{j=N}a_j.y^{[j/m)]+\epsilon(j,m)}.(dy/y)^{\otimes j}\otimes dx^{\otimes (N-j)},$$ avec les $a_j$ holomorphes, o\`u $[*]$ est la partie enti\`ere, et o\`u $\epsilon(j,m):=0$( resp. $1$) si $m$ divise $j$ (resp. sinon).

Notons que $[j/m]+\epsilon(j,m)$ vaut $0$ si $j=0$; vaut $1$ si $1\leq j\leq m$; et vaut $k$ si et seulement si $(k-1)m+1\leq j\leq km$.

\end{definition}

\

Afin de faciliter l'\'ecriture, on pose en 4.2-5 ci-dessous: $[Sym^N](\Omega_{(B/\Delta)}^1):=S^N(m,B_F)=S^N(m)$, donc $S^N(\infty)$:=$Sym^N(\Omega_B^1(log D))$.

\

\begin{remark} {\it 

\

$\bullet$ L'expression ci-dessus conserve la m\^eme forme si l'on change les coordonn\'ees locales $(x,y)$, avec $y=0$ \'equation locale de $D$.

$\bullet$ Rappelons que $\Omega_B^1(log D)$ est le faisceau pr\'ec\'edent, lorsque $m=\infty$, et que $[Sym^N](\Omega_{(B/\Delta)}^1)\subset Sym^N(\Omega^1_{B}(log D))$. 

\

Plus g\'en\'eralement, si $\Delta=(1-1/m).D$ et si $\Delta'=(1-1/m').D$, pour $1\leq m\leq m'\leq \infty$, alors $[Sym^N](\Omega_{(B/\Delta)}^1)\subset [Sym^N](\Omega_{(B/\Delta')}^1))$.}

$\bullet$ L'application naturelle de multiplication: $S^N(m)\otimes S^M(m)\to S^{N+M}(m)$, qui envoit :

$[(\sum_{j=0}^{j=N}a_j.(dy^{\otimes j}/y^{[(1-1/m)j]})\otimes dx^{\otimes (N-j)})]\otimes [\sum_{k=0}^{k=M}b_k.(dy^{\otimes k}/y^{[(1-1/m)k]})\otimes dx^{\otimes (M-j)}]$ sur: $[\sum_{j=0}^{j=N+M}c_h.(dy^{\otimes h}/y^{[(1-1/m)h]})\otimes dx^{\otimes (N+M-h)}]$, avec $[c_h:=\sum_{j+k=h}a_j.b_k.y^{r(j,k;m)}]$  est bien d\'efinie, puisque:

 $r(j,k;m):=[j/m]+\epsilon(j,m)+[k/m]+\epsilon(k,m)-\lbrace[(j+k)/m]+\epsilon(j+k,m)\rbrace\geq 0$, pour tous $j,k,m$, par un calcul imm\'ediat.
\end{remark}

\

\begin{proposition} Si $\cal O$$_D(D)\cong \cal O$$_D$, et si $g$ est le genre de $D$, alors:

$h^0(B,S^{qm}(m))\geq h^0(B,S^{qm}(\infty))-((2(g-1)).m^2q^3)/6).-A(m)q^2$, o\`u $A(m)$ est une constante ne d\'ependant que de $m$.

\end{proposition}

{\bf D\'emonstration:} Le faisceau quotient $Q:=S^N(\infty)/S^N(m)$ est support\'e par $D$, et 
$h^0(B,S^{qm}(m))\geq h^0(B,S^{qm}(\infty))-h^0(D,Q).$

Nous allons montrer que $h^0(D,Q)\leq ((2(g-1)).m^2q^3)/6).-A(m)q^2,$ ce qui \'etablira l'assertion.

Le faisceau $Q$ admet une filtration naturelle $Q_N\subset Q_{N-1}\subset ....\subset Q_1\subset Q_0=Q$ telle que, pour $j=0,1,...,N$, on ait: $ Q_{(j-1)}/Q_j=Sym^{N-j}(\Omega_D^1)\otimes F_j,$ o\`u $F_j$ est un faisceau localement libre de rang $r_j:=[j/m]+\epsilon(j,m)$, et admettant une filtration $F_{j,0}\subset F_{j,1}\subset ...\subset F_{j,r_j}=F_j$, avec $F_{j,h}/F_{j,h-1}\equiv (N^*)^{\otimes h}$, o\`u $N^*$ est le dual du fibr\'e normal \`a $D$ dans $B$.

Nous n'utiliserons ces filtrations que lorsque $N\cong \cal O$$_D$.

Dans ce cas, on a donc: 

$h^0(D,Q_j)\leq r_j.h^0(D,K_D^{\otimes(N-j)})\leq  (2(N-j)-1)(g-1).r_j$, si $N-j\geq 2$.

\

On suppose d\'esormais que $N=qm$ est un multiple de $m$. 

\

On a alors (puisque $1\leq h:=r_j\leq q$ si $j=(h-1).m+k,$ avec $1\leq k\leq m$): 

$$h^0(D,Q)\leq \sum_{j=1}^{j=N}h^0(D,Q_j)=
1+ \sum_{h=1}^{h=q}h.[ \sum_{k=1}^{k=m}h^0(D,K_D^{\otimes ((q-h+1)m-k)})]$$
$$\leq \sum_{h=1}^{h=q}mh.(2(q-h+1)m-1)(g-1)$$

Cette derni\`ere somme vaut: $(2(g-1)m^2.q^3)((1/2)-(1/3))-A(m)q^2$, o\`u $A(m)$ d\'epend de $m$ seulement, ce qui \'etablit l'assertion.

\

\begin{proposition} Soit $B,D,m$ comme ci-dessus, et $F:=\Omega_B^1(log D)$. 

On pose: $e_1:=c_1(F)$ et $e_2:=c_2(F)$, ainsi que $c_1:=c_1(B)$ et $c_2:=c_2(B)$. 

On suppose que $\cal O$$_D(D)\cong \cal O$$_D. ($ Donc $2(g-1)=K_B.D)$. 

Alors:

1. $e_1\equiv -c_1+D$, $e_1^2=c_1^2+2K_B.D$. 

2. $e_2=c_2+K_B.D$, $e_1^2-e_2=c_1^2-c_2+K_B.D$

3. $h^0(B,S^{qm}(m))\geq \alpha. (mq)^3/6-O(q^2)$, si  $\alpha:=(c_1^2-c_2)+(1-1/m).K_B.D>0$, et si $\kappa(B)=1$ et $g(D)\geq 2$.

\end{proposition}

{\bf D\'emonstration:} Le faisceau quotient $R:=\Omega_B^1(log D)/\Omega_B^1$ est support\'e sur $D$  et isomorphe \`a $\cal O$$_D$. Donc $e_1=-c_1+D$ (restreindre \`a une courbe tr\`es ample arbitraire).

Pour calculer $e_2$, on applique la formule de Riemann-Roch \`a l'\'egalit\'e $\chi(\Omega_B^1(log D))=\chi(\Omega_B^1)+\chi(\cal O$$_D)$, en observant que $K_B.D=-2\chi(\cal O$$_D$). 

\

On obtient ainsi la seconde assertion, puisque cette \'egalit\'e s'\'ecrit:

$(e_1^2/2)-e_2-K_B.e_1/2=((c_1^2/2)-c_2+K_B.c_1/2)-K_B.D/2,$ 

que $e_1^2/2=c_1^2/2+K_B.D,$ et que $-K_B.e_1/2=K_B.c_1/2-K_B.D/2$.

\

La derni\`ere assertion r\'esulte des deux propositions pr\'ec\'edentes, de l'\'egalit\'e de Riemann-Roch: $\chi(B,Sym^N(\Omega_B^1(log D)))=((e_1^2-e_2)/6).N^3+O(N^2)$, et du lemme ci-dessous, en remarquant que:

 $h^0(B,(Sym^M(\Omega^1_B(log D))))\geq \chi (B,(Sym^M(\Omega^1_B(log D))))-(1+g(D))$, puisque:
 
 $h^2(B,(Sym^M(\Omega^1_B(log D))))\leq 1+g(D)$.
 
 \

\begin{lemma} Soit $B$ une surface projective complexe lisse et connexe, et $D\subset B$ une courbe projective lisse et connexe de $B$. On suppose que $g(D)\geq 2$, et que $N:=\cal O$$_D(D)\cong \cal O$$_D$. Soit $T_B(-log D)$ (resp. $T_B(-D))$ le faisceau des germes de champs de vecteurs holomorphes sur $B$ qui sont tangents \`a $D$ (resp. nuls le long de $D$).

Alors, pour tout $M>0$:

1. L'injection naturelle suivante a un conoyau de dimension $1+g(D)$ au plus:

$H^0(B,(Sym^M(T_B(- D)))\otimes K_B) \to H^0(B,(Sym^M(T_B(-log D)))\otimes K_B)$.

2. Si $B$ est une surface elliptique, alors:

$h^2(B,(Sym^M(\Omega^1_B(log D))))=h^0(B,(Sym^M(T_B(-log D)))\otimes K_B)\leq 1+g(D)$.

\end{lemma}

\

{\bf D\'emonstration:} Soit $T:=T_B(-log D)$, et $T_{\vert D}$ sa restriction \`a $D$ (comme faisceau localement libre de rang $2$).

\

On a deux suites exactes naturelles:

$0\to T\to T_B\to N\to 0$, et par restriction \`a $D$: 

$0\to T_D\to T_{\vert D}\to N(-D)\cong \cal O$$_D\to 0$.

\

Donc $Sym^M(T_{\vert D})\otimes K_B$ admet une filtration croissante par les $F_j,j=0,1,...,M$ telle que $F_j/F_{j-1}\cong (T_D)^{\otimes j}\otimes K_D$, de sorte que $H^0(D,F_j/F_{j-1})=0$ pour $j\leq M-2$. 

Et que $H^0(D,F_j/F_{j-1})=1$ (resp. $g(D))$ si $j=M-1$ (resp. $j=M$). 

Donc: $h^0(D,Sym^M(T_{\vert D})\otimes K_D)\leq \sum_{0\leq j\leq M}h^0(D,F_j/F_{j-1})=1+g(D)$. 

D'o\`u la premi\`ere assertion, en consid\'erant la suite exacte: 

$H^0(B,(Sym^M(T_B(- D)))\otimes K_B) \to H^0(B,(Sym^M(T_B(-log D)))\otimes K_B)\to  H^0(D,(Sym^M(T_B(-log D)_{\vert D}))\otimes K_D)$.

\

Pour la seconde assertion, on note $g:B\to G$ la fibration elliptique de $B$. Alors $g:D\to G$ est finie, et $H^0(F,(Sym^M(T_B(-log D)_{\vert F}))\otimes K_{B\vert F})=0$, si $F$ est une fibre g\'en\'erique (lisse elliptique) de $g$, puisque $K_{B\vert F}\cong \cal O$$_F$.

Donc $g_*(Sym^M(T_B(- D)))\otimes K_B)=0=H^0(B,(Sym^M(T_B(-D)))\otimes K_B)$, et enfin: $h^0(B,(Sym^M(T_B(-log D)))\otimes K_B)\leq h^0(B,(Sym^M(T_B(- D)))\otimes K_B)+(1+g(D)\leq 1+g(D)$.

\

 Le corollaire 4.4  r\'esulte imm\'ediatement des assertions  pr\'ec\'edentes.

\

\subsection{Un th\'eor\`eme de confinement pour les courbes holomorphes $m$--tangentes \`a $D$}

\

Soit $B_0,D_0$ une surface munie d'une courbe $D_0$, poss\'edant les propri\'et\'es \'enonc\'ees dans le \S \ref{bt}. On \'eclate $B_0$ en $F=D_0\cap D'_0$, pour un choix g\'en\'erique de $D'_0$, obtenant $B,D$.

\noindent Soit $h:\bC\to B$ une courbe enti\`ere $m$-tangente \`a $D$, ou encore qui d\'efinit un morphisme orbifolde vers $(B, \Delta)$, c'est-\`a-dire telle que
$$h^{*}D= \sum_j\nu_j\delta_{z_j}$$
o\`u toutes les multiplicit\'es $\nu_j$ sont {\it divisibles} par $m$. (Voir cependant la remarque 4.10).

Nous allons montrer que (pour un choix g\'en\'erique de $D'_0)$, les points d'intersection de $h(\bC)$ avec $D$ se projettent en un ensemble fini $G:=G(F)$ de $D_0$, ind\'ependant, non seulement de $h$, mais aussi de $F$. Ces faits \'etablis, on peut conclure exactement comme dans le cas de la version alg\'ebrique de l'hyperbolicit\'e. Nous ne reviendrons donc pas sur cet argument.

\

La finitude et l'ind\'ependance en $h$ seront d\'emontr\'ees dans un premi\`ere \'etape. L'ind\'ependance en $F$ sera \'etablie dans une seconde \'etape, qui n\'ecessite de raffiner les arguments de la premi\`ere \'etape, et fait appel \`a des propri\'et\'es plus fines des courants, rappel\'ees pour la commodit\'e du lecteur.

\

{\bf Premi\`ere \'etape: Construction de $G(F)$} 

\

Consid\'erons \'egalement  un fibr\'e ample $A$ sur la surface 
$B$, muni d'une m\'etrique 
$\psi$ \`a courbure positive. D'apr\`es le corollaire 4.6, il existe une
section 
$\omega\in H^0(B, [Sym^N](\Omega ^1_{(B/\Delta)})\otimes A^{-1})$, non-nulle.

On utilise \'egalement le fibr\'e
$A$ pour d\'efinir la fonction caract\'eristique $r\to T_A(h, r)$ de la courbe $h$ comme suit
(voir \cite{[G]} pour une pr\'esentation plus d\'etaill\'ee de cette notion) 

$$T_A(h, r):= \int_0^r{{dt}\over {t}}\int_{D(t)}h^*\Theta_\psi(A)$$
pour chaque $r\in \bR_+$.

L'hypoth\`ese sur le diviseur $h^*D$ ci-dessus
montre qu'on peut faire agir $\omega$ sur la d\'eriv\'ee de la courbe (ou application) $h$, obtenant une section holomorphe de $h^*A^{-1}$. 
Le lemme suivant montre  
que $\omega$ fournit une \'equation diff\'erentielle alg\'ebrique pour la courbe $h$.

\begin{lemma} La section $\omega(h')$ 
du fibr\'e $h^*A^{-1}$ s'annule identiquement sur $\bC$.

\end{lemma}

\noindent {\bf Preuve:} La d\'emonstration r\'epose sur la proposition suivante,
due \`a P. M Wong et Y.T Siu (voir \cite{[W]}, \cite{[R]}, \cite{[S]}). 
On remarquera que les arguments classiques de courbure 
n\'egative ne fonctionnent dans ce cadre qu'en utilisant une version 
orbifolde ad\'equate de la notion de m\'etrique hermitienne 
sur un fibr\'e en droites (approche que nous n'avons pas d\'evelopp\'ee ici), 
ceci \`a cause des p\^oles de $\omega$.

\begin{proposition} Soit $\omega$ une diff\'erentielle 
m\'eromorphe sur une vari\'et\'e projective $B$, 
qui admet des p\^oles logarithmiques le long d'un diviseur \`a croisements normaux $D$. 
Soit $h:\bC\to B$ une courbe enti\`ere d'image non contenue dans $D$.
Alors on a:
$$\int_0^{2\pi}\log^+\vert \omega(h'(r\exp(i\theta)))\vert _{\psi^{-1}} d\theta\leq C\log(T(h, r))$$
pour tout r\'eel positif $r\in \bR\setminus E$, o\`u $E$ est un ensemble de mesure de Lebesgue finie.  
\end{proposition}

\

Afin de donner une id\'ee de la d\'emonstration de 4.8, nous rappelons le classique lemme de la d\'eriv\'ee logarithmique de Nevanlinna, voir \cite{[N]}.

\begin{theorem} Soit $\tau:\bC\to\bP^1$ une application m\'eromorphe. Alors il existe 
un ensemble $E\subset \bR_+$ de mesure de Lebesgue finie, tel que 
pour tout $r\in \bR\setminus E$ on ait
$$\int_0^{2\pi}\log^+\Big\vert 
{{\tau'(r\exp(i\theta))}\over {\tau(r\exp(i\theta))}}\Big\vert d\theta\leq C\log(T(h, r)).$$
\end{theorem}  

\

La proposition 4.8 est 
une cons\'equence assez immediate 
du th\'eor\`eme pr\'ec\'edent.
En effet, on peut supposer le fibr\'e $D$ tr\`es ample,
quitte \`a lui ajouter un multiple assez grand de $A$.
Il existe donc des sections $(\sigma_j)_{j=1,...N}$ de $D$ telles qu'au voisinage de tout point
$x\in B$, on puisse extraire de l'ensemble $\displaystyle 
\log\Bigl({{\sigma_j}\over{\sigma _k}}\Bigr)_{j, k= 1,...,N}$
des coordonn\'ees locales holomorphes.

Consid\'erons maintenant une courbe enti\`ere
$h_1:\bC\to \bP^1$. Le th\'eor\`eme de Nevanlinna peut 
\^etre r\'eformul\'e comme suit: la moyenne sur les cercles
du logarithme de $\displaystyle \vert h_1^*{{dz}\over {z}}/dt\vert$
est domin\'ee par ${\cal O}(\log T(h_1, r))$. Enfin,
la fonction $t\to \vert \omega(h'(t)))\vert _{\psi^{-1}}$ s'exprime localement
sous la forme: ${\cal P}\Bigl(h(t), {{d}\over {dt}}(\log h_{jk}(t))\Bigr)$,
o\`u ${\cal P}$ est un polyn\^ome, 
et $\displaystyle h_{jk}:= {{\sigma_j}\over {\sigma_k}}\circ h$. La compacit\'e de $X$ et le
lemme de la d\'eriv\'ee logarithmique impliquent alors le r\'esultat (4.8); pour plus de d\'etails, 
voir \cite{[R]}, pp. 289-302.

\

Maintenant, si par l'absurde $\omega(h')$ n'est pas identiquement nulle, on a: 
$$i\partial\overline\partial \log\vert \omega(h')\vert ^2_{\psi^{-1}}\geq h^*\Theta_{\psi}(A)$$
par l'\'equation de Poincar\'e-Lelong. Si on int\`egre l'in\'egalit\'e ci-dessus \`a la mani\`ere de 
Nevanlinna, on obtient
$$\int_0^r{{dt}\over{t}}\int_{D(t)}i\partial\overline\partial \log\vert \omega(h'))
\vert ^2_{\psi^{-1}}\geq
T_A(h, r)$$
pour tout r\'eel positif $r$. La formule de Jensen donne ensuite
$$\int_0^{2\pi}\log^+\vert \omega(h'(r\exp(i\theta)))\vert _{\psi^{-1}} d\theta
\geq T_A(h, r)+ {\cal O}(1),$$
et compte tenu du fait que notre op\'erateur $\omega$ est 
(en particulier) \`a p\^oles logarithmiques le long de $D$,
le r\'esultat de Wong combin\'e avec l'inegalit\'e pr\'ecedente donnent
$$\log(T_A(h, r))\geq T_A(h, r)+ {\cal O}(1) $$
pour tout $r\in \bR\setminus E$. Lorsque $r\to \infty$, on obtient une contradiction, et le lemme 4.7 est ainsi d\'emontr\'e.

\

On se place maintenant \`a nouveau dans la situation $(B,D,N,F,D_0,m)$ et les hypoth\`eses du corollaire 4.6.

\

\begin{proposition}\label{fin}

(1) Il existe un sous-ensemble fini $G=G(F)$ de points $(p_j)_{j\in J}$ de la courbe $D=D_0$, tels que pour toute courbe enti\`ere $h:\bC\to B$ $m$-tangente \`a $D$, l'intersection (ensembliste) de l'image de $h$ avec $D$ est contenue dans la r\'eunion finie $G=G(F)$ des $(p_j)$. 

(2) L'ensemble $G(F)$ est ind\'ependant de $F$.

\end{proposition}

\

Ceci acquis, la preuve sera finie, car on peut reprendre
les arguments utilis\'es dans le cadre alg\'ebrique: la construction pr\'esent\'ee dans le paragraphe
3 montre qu'il existe une surface de type g\'en\'eral $S'$, telle que toute courbe (alg\'ebrique ou non) 
dans $B$ , $m$-tangente \`a $D$, parametr\'ee par $\bC$, et  dont l'image ne contient aucun des points $(p_j)$ se rel\`eve \`a $S'$. 

Mais l'in\'egalit\'e $(c_1^2-c_2)>0$ sur les classes de Chern de $S'$, et le th\'eor\`eme de McQuillan (\cite{[MQ]}) montrent que
l'image de $h$ est contenue dans une courbe alg\'ebrique. Les r\'esultats du
paragraphe pr\'ecedent ach\`event alors la d\'emonstration, puisque la compos\'ee $\varphi\circ h:\bC\to X$ est $m$-tangente \`a $D$, pour toute application holomorphe $h:\bC\to X$.

\

{\bf D\'emonstration:} Pour montrer l'existence des points $(p_j)$ dans \ref{fin} (1), observons tout d'abord que 
$\omega$ peut \^etre vue comme une section m\'eromorphe de ${\cal O}_P(mq)\otimes A^{-1}$, o\`u $P:=\bP(T_B)$.

Supposons que l'image de la courbe $h$ ne soit pas contenue dans $D$. Alors le lemme 4.7
montre que l'image de la courbe d\'eriv\'ee $h':\bC\to \bP(T_B)$ est contenue
dans une composante irr\'eductible $Y$ de l'ensemble des z\'eros de $\omega$. 
On analyse les deux cas qui peuvent apparaitre.

\

1. La dimension de la projection $\pi(Y)$ sur $X$ est 1; alors en particulier,
l'image de $h$ est contenue dans une courbe alg\'ebrique, cas trait\'e auparavant.

\

2. Si, enfin, $\pi_{_\vert Y}:Y\to B$ est une application (g\'en\'eriquement) finie, il existe un nombre fini de points de $B$ pour lesquels la fibre de  $\pi_{_\vert Y}$ est de dimension positive.
De plus, sur $Y$ (ou plut\^ot, sa d\'esingularis\'ee), on a un feuilletage canonique ${\cal F}$, 
dont les disques tangents sont d\'ecrits comme suit.

Soit $\tau:(\bC, 0)\to B$ un germe de disque holomorphe, dont l'image de la
d\'eriv\'ee projectivis\'ee $\tau'$ est contenue dans $Y$ 
(autrement dit, $\tau$ est tangente au multi-feuilletage
d\'efinit par $Y$). Alors la d\'eriv\'ee de $\tau'$ \`a l'origine est 
tangente \`a ${\cal F}$. 

Il existe encore deux sous-possibilit\'es. On note $\hat D$ la r\'elev\'ee canonique de 
la courbe $D$ dans $\bP(T_Y)$. Si $\hat D$ est contenue dans $Y$, c'est une courbe invariante de 
${\cal F}$, et donc $h'$ ne peut intersecter $\hat D$ qu'aux points singuliers de ${\cal F}$ (sinon, $h'$ est contenue dans $\hat D$, par le th\'eor\`eme d'unicit\'e
des solutions pour les \'equations diff\'erentielles ordinaires). 

Dans l'autre cas, $\hat D$ intersecte 
$Y$ en un nombre fini des points.

Dans les deux cas, les points $(p_j)\subset D$ 
sont obtenus par projection de $P$ sur $B$, soit des singularit\'es de ${\cal F}$, soit de
l'intersection de la relev\'ee canonique de $D$ avec $Y$. Ce qui ach\`eve la preuve de \ref{fin} (1).
\

\begin{remark} La preuve du lemme 4.7 fonctionne sans changement notable dans le cadre plus 
naturel des courbes enti\`eres $m$-tangentes \`a $D$ au sens {\it non-classique}. En effet, il suffit d'observer que 
dans ce cas $\omega(h')$ est bien d\'efinie seulement en dehors du support du diviseur
$\sum \nu_j\delta_{z_j}$ (on emploie ici les notations du d\'ebut de cette section), 
et de plus la norme de cette section est born\'ee. On invoque ensuite
le th\'eor\`eme de prolongement des fonctions psh, et le reste des arguments est identique.
\end{remark}

\

{\bf Seconde \'etape: $G(F)$ ne d\'epend pas de $F$.}

\

{\bf D\'emonstration (suite):} On va \'etablir l'ind\'ependance des points
$(p_j)$ par rapport \`a la configuration $F$ utilis\'ee pour construire par \'eclatements la surface 
$B_F$, et \'etablir ainsi \ref{fin} (2).

La d\'emonstration va aussi se faire en deux \'etapes: dans la premi\`ere (\ref{mult}), on montre que la conclusion est v\'erifi\'ee si la multiplicit\'e du courant d'Ahlfors sur $B$ associ\'e \`a $h$ a une multiplicit\'e nulle aux points d'intersection de $h(\bC)$ avec $D$. Dans la seconde (\ref{kod}), on \'etablit cette propri\'et\'e comme cons\'equence du fait que la dimension de Kodaira de $B$ est $1$.

\

Pour ceci, on a besoin de faire quelque rappels concernant les courants
associ\'es aux courbes enti\`eres, et leurs propri\'et\'es num\'eriques. 
Les r\'ef\'erences principales qu'on utilisera sont [MQ98], [Br00] et [Pa03].

Soit $X$ une vari\'et\'e projective, et soit $A\to X$ un fibr\'e ample, muni
d'une m\'etrique $\psi$ dont la forme de courbure $\Theta _\psi(A)$ est d\'efinie positive sur 
$X$. Si $h:\bC \to X$ est une courbe enti\`ere, alors
l'indicatrice de croissance de 
$h$ par rapport \`a $\Theta_\psi(A) $ est d\'efinie par
$$ T_A(h, r):=  
\int_0^r{{dt}\over {t}}\int_{\Delta(t)}h^*\Theta_\psi(A)$$
o\`u $\Delta(t)\subset \bC$ est le disque de rayon $t$.
McQuillan associe \`a $h$ un courant positif 
ferm\'e de la fa\c{c}on suivante: pour tout r\'eel
positif $r$, consid\'erons le courant de type $(n- 1, n-1)$ d\'efini par
$$(T_r[h], \alpha):= {{\int_0^r{{dt}\over {t}}\int_{\Delta(t)}h^*\alpha}\over
{T(h, r)}},  
$$ 
o\`u $\alpha$ est une forme de type $(1,1)$, de classe ${\cal C}^\infty$ sur $X$.
Maintenant, si $(r_k)$ est une suite qui tend vers l'infini, la famille de courants 
$\displaystyle T_{r_k}[h]$ est de masse uniform\'ement born\'ee, donc par pr\'ecompacit\'e,
admet une limite faible, not\'ee $T[h]$. C'est un courant positif, 
et de plus, quitte \`a bien choisir la suite $(r_k)$, il sera ferm\'e (cette affirmation est
une cons\'equence des r\'esultats classiques d'Ahlfors et Nevanlinna).
On rappelle tr\`es bri\`evement les propri\'et\'es num\'eriques de ces 
courants qui vont servir par la suite.

\

1. Si la courbe $h$ est Zariski dense, et si $\{\alpha \}$ 
est une classe contenant un courant positif ferm\'e, alors
$(T[h], \alpha)\geq 0$ (voir [MQ98], [Pa03]).

\

2. On d\'efinit {\sl la multiplicit\'e asymptotique} $mult (T[h], x)$ de la courbe $h$ au point $x\in X$ comme suit:
soit $\hat X$ l'\'eclatement de $X$ en $x$, et soit  $\hat h$ la relev\'ee de la courbe $h$.
Si $E$ d\'esigne le diviseur exceptionnel, alors la multiplicit\'e de 
$h$ en $x$ est par d\'efinition $mult (T[h], x):= (T[\hat h], c_1(E))$. 
Le r\'esultat dont dont nous aurons  besoin est le suivant:
{\sl si $h$ est une courbe enti\`ere Zariski dense
trac\'ee sur une surface projective, et dont la multiplicit\'e asymptotique en un point
est positive, alors la classe de cohomologie du courant 
$T[h]$ contient un courant k\"ahl\'erien} (voir [MQ05] et \'egalement [Pa03]).
(Le r\'esultat analogue en g\'eom\'etrie alg\'ebrique est le suivant:
soit $C$ une courbe nef sur une surface projective $S$, et soit $x\in C$; on suppose que
la transform\'ee propre $\hat C$ de $C$ par rapport \`a l'\'eclatement de 
$S$ en $x$ est d'autointersection positive. Alors la courbe $C$ est ``big" (lorsque vue comme diviseur sur $S)$.

\

Apr\`es ce bref d\'etour, revenons aux ensembles $F\subset B_0$, et aux surfaces 
$B_F$ correspondantes. On fixe une de ces surfaces $B=B_{F_0}$ associ\'ee \`a $F_0$ comme surface 
de r\'ef\'erence, et \`a toute courbe $h_{F}:\bC \to B_{F}$, on associe la courbe 
$h_{F,0}:\bC \to B$ en projettant $h_F$ sur $B_0$, et en consid\'erant la 
transform\'ee propre de cette projection sur $B_{F_0}$. On remarquera que si 
la courbe $h_F$ est $m$-tangente \`a $D_F$, alors 
la courbe $h_{F,0}$ est $m$--tangente \`a $\displaystyle D_{F_0}$, sauf (\' eventuellement) 
en un nombre fini de points
contenus dans $F\cup F_0$. On utilise cette observation pour montrer le lemme suivant
(qui renforce 4.10). Dans l'enonc\'e suivant, on emploie les 
notations et hypoth\`eses du debut du paragraphe 4.4.

\begin{lemma}\label{mult} Soit $F\subset D_0$ une configuration, et soit 
$h_F:\bC \to B_F$ une courbe enti\`ere Zariski dense, telle
que $mult(T[h_{F,0}], x)= 0$, en tout point
$x\in B$. Alors
la section m\'eromorphe $\omega(h_{F,0}')$ 
du fibr\'e $h_{F,0}^*A^{-1}$ s'annule identiquement sur $\bC$.

\end{lemma} 

\noindent {\bf Preuve:} La d\'emonstration est identique \`a celle 
du lemme 4.10; la non-annulation de la section m\'eromorphe $\omega(h_{F,0}')$ et 
l'\'equation de Poincar\'e--Lelong montrent que
$$i\partial\overline\partial \log\vert \omega(h_{F,0}')\vert ^2_{\psi^{-1}}
+N\sum _p \sum _{j} \nu_{j,p}\delta_{z_{j, p}}
\geq h_{F,0}^*\Theta_{\psi}(A)$$
o\`u l'on a not\'e $(s_p)_{p=1..M}$ l'ensemble fini des points du diviseur $D$ ou 
$h_{F,0}$ n'est pas $m$--tangente \`a $D$,
et les quantit\'es intervenant dans la seconde somme sont d\'efinies par 
$\displaystyle h_{F,0}^{-1}s_p= \sum_j \nu_{j,p}z_{j,p}$  
. En int\'egrant cette in\'egalit\'e, 
on obtient:

$$\int_0^r{{dt}\over{t}}\int_{D(t)}i\partial\overline\partial \log\vert \omega(h_{F,0}'))
\vert ^2_{\psi^{-1}}+ N\sum_{p=1}^M\sum_{0<\vert z_{j,p}\vert < r}\nu_{j,p}\log{{r}\over {\vert z_{j,p}\vert }}\geq
T_A(h_{F,0}, r)$$
pour tout r\'eel positif $r$. La formule de Jensen donne ensuite
$$\int_0^{2\pi}\log^+\vert \omega(h'(r\exp(i\theta)))\vert _{\psi^{-1}} d\theta
+ N\sum_p\sum_{0<\vert z_{j,p}\vert < r}\nu_{j,p}\log{{r}\over {\vert z_{j,p}\vert }}
\geq T_A(h_{F,0}, r)+ {\cal O}(1),$$
On divise l'inegalit\'e pr\'ec\'edente par l'indicatrice de croissance
$T_A(h_{F,0}, r)$, et on fait $r\to \infty$; on en d\'eduit

$$1 < N\sum_{p=1}^M mult(T[h_{F,0}], s_p).$$ 

Mais par hypoth\`ese, la multiplicit\'e de la courbe en un point arbitraire 
vaut z\'ero, et la contradiction ainsi obtenue d\'emontre 
la proposition.

\

On remarquera ici que la forme $\omega$ 
est ind\'ependente de $F$.

\

L'\'etape suivante de la preuve est de montrer qu'on a 
$mult(T[h_{F,0}], x)= 0$, pour tout $x\in B$, si $h_{F,0}$ est Zariski dense.

C'est une cons\'equence du r\'esultat g\'en\'eral suivant.

\begin{lemma}\label{kod} Soit $B$ une surface projective non--singuli\`ere,
et soit $h:\bC\to B$ une courbe enti\`ere Zariski dense, tangente \`a un
multi-feuilletage. Si $mult(T[h], x)> 0$ pour un point $x\in B$, alors 
$\kappa(B)\leq 0$.

\end{lemma}

\noindent {\bf Preuve} Cet enonc\'e est implicite 
dans les articles [MQ98], [Br00], mais pour la commodit\'e du lecteur, on 
va en donner la d\'emonstration.

Par hypoth\`ese, l'image de la 
courbe d\'eriv\'ee $h'$ de $h$ est contenue dans un 
diviseur $Y\subset \bP(T_B)$, et $h'$ est tangente \`a un feuilletage sur $Y$.

Quitte \`a modifier $Y$, on peut supposer qu'on se trouve dans 
la situation suivante: il existe $\pi: B_1\to B$ une application surjective,
(compos\'ee d'une application finie avec une suite d'\'eclatements),
et une courbe $h_1:\bC\to B_1$, Zariski dense, tangente \`a un feuilletage 
r\'eduit au sens de Seidenberg (voir e.g. [MQ98] et les r\'ef\'erences dans cet article), et telle que
$\pi\circ h_1= h$.
Dans ce cas, les r\'esultats de McQuillan et Brunella ([MQ98], [Br00]) montrent que
$$\displaystyle (T[h_1], K_{B_1})\leq 0$$ 
Par ailleurs, on a $\displaystyle K_{B_1}= \pi^*K_B+ R$, o\`u
$R$ est un diviseur effectif, et la propri\'et\'e (1) cit\'ee auparavant montre que 
$(T[h_1], c_1(R))\geq 0$; on en d\'eduit ainsi l'in\'egalit\'e:
 
$$\displaystyle (T[h_1], \pi^*K_{B})\leq 0$$ 

Supposons \`a pr\'esent que $\kappa(B)\geq 1$. Alors le diviseur canonique
admet une d\'ecomposition de Zariski $K_B= P+E$, o\`u $P$ est un diviseur nef, 
de dimension num\'erique au moins 1, et $E$ est un diviseur effectif.
La relation pr\'ec\'edente montre que  $\displaystyle (T[h], P)\leq 0$
(car l'image directe du courant $T[h_1]$ par $\pi$ est simplement $T[h]$).
Par hypoth\`ese, la courbe $h$ admet une multiplicit\'e positive en 
$x$, et la propri\'et\'e 2 montre que $\displaystyle (T[h], P)> 0$.
Le lemme est d\'emontr\'e.

\

En conclusion des r\'esultats pr\'ec\'edents, 
il existe une section $\omega\in H^0(B, [Sym^N](\Omega ^1_{(B/\Delta)})\otimes A^{-1})$
telle que pour tout ensemble $F$, et toute courbe enti\`ere 
$h_F:\bC\to B_F$ Zariski dense et $m$-tangente \`a $D_F$, on ait $\omega (h'_{F,0})= 0$.

Donc il existe un ensemble fini $G$ de points $(p_j)$ tel que pour tout ensemble $F$,
l'intersection de la projection des courbes enti\`eres Zariski denses de $X_F$
sur $B_F$ est contenue dans la r\'eunion des $(p_j)$. La proposition \ref{fin}, et donc aussi le th\'eor\`eme 4.3 sont ainsi compl\`etement d\'emontr\'es.

\vspace{1cm}
\small
\begin{tabular}{lcl}
Fr\'ed\'eric Campana && Mihai P\u aun \\ D\'epartement de
Math\'ematiques&& D\'epartement de
Math\'ematiques \\
Universit\'e de Nancy && Universit\'e de Nancy\\
F-54506 Vandoeuvre-les-Nancy, France &&  F-54506 Vandoeuvre-les-Nancy, France \\
campana@iecn.u-nancy.fr && mihai.paun@iecn.u-nancy.fr\\
\end{tabular}

\end{document}